\documentclass[reqno,fleqn]{amsart}
\usepackage{amsfonts,amsthm,amssymb}
\usepackage{lipsum}
\usepackage[utf8]{inputenc}
\usepackage[english]{babel}

\usepackage{lmodern,textcomp}
\usepackage{amsmath}
\usepackage{mathtools}
\usepackage{latexsym}
\usepackage{tikz}
\usepackage{fancyvrb}
\usepackage{epsfig}
\usepackage{pstricks,slashbox,multirow}
\usepackage{graphicx}
\usepackage{rotating}
\usetikzlibrary{positioning}
\usepackage{subcaption}
\usepackage{datetime}
\newtheorem{theorem}{Theorem}[section]

\newtheorem{definition}[theorem]{Definition}

\newtheorem{rem}[theorem]{Remark}

\title[On Type-2 isomorphic $C_n(R)$. PART 9: Computer programs to show Type-2 isomorphic $C_n(R)$]{A study on Type-2 isomorphic circulant graphs. \\
PART 9: Computer programs to show Type-1 $\&$ -2 isomorphic circulant graphs}
\author{\sc Vilfred~ Kamalappan} 
\address{Department ~of~ Mathematics, ~Central ~University ~of~ Kerala,~Periye,  ~Kasaragod,~ Kerala, \linebreak India~ - ~671 316.}
\email{vilfredkamal@gmail.com}
\author{\sc Wilson~ Peraprakash} 
\address{Department~ of~ Mathematics,~ S.T.~ Hindu~ College,~Nagercoil,  ~Kanyakumari~ District,~ Tamil ~Nadu, \linebreak India~-~629 002.}
\email{wilsonperapras@gmail.com}

\date{}

\begin{document}

\begin{abstract} Elspas and Turner \cite{eltu} raised a question on the isomorphism of $C_{16}(1,3,7)$ and $C_{16}(2,3,5)$ and Vilfred \cite{v96} gave its answer by defining Type-2 isomorphism of $C_n(R)$ w.r.t. $m$ $\ni$ $m$ = $\gcd(n, r) > 1$ and also modified to $m$ and $m^3$ divides $\gcd(n, r)$ and $n$, respectively, $r\in R$ and $r,n\in\mathbb{N}$ and studied such graphs \cite{v13,v20, v2-1}. Using a $C^{++}$ computer program, the authors obtained families of Type-2 isomorphic $C_{n}(R)$ w.r.t. $m$ = 2,3,5,7 for  $n\in\mathbb{N}$ as well as $C_{np^3}(R)$ w.r.t. $m$ = $p$ where $p$ is a prime and $n\in\mathbb{N}$. In this paper, we present the $C^{++}$ program and also a VB program POLY415.EXE which is used to show how Type-1 and Type-2 isomorphisms of a circulant graph take place as well as for checking and finding Type-1 and Type-2 circulant graphs of a given order. This article is the $9^{th}$ part of a detailed study on Type-2 isomorphic circulant graphs having ten parts \cite{v2-1}-\cite{v2-10}. 
\end{abstract}

\maketitle

\section{Introduction}

Graphs are powerful data structure to represent objects and their concepts. Objects are nothing but vertices and edges describe relation among objects. Deciding whether two graphs are structurally identical or isomorphic, is a classical algorithmic problem \cite{mp20} that has been studied since the early days of computing. When the graphs are not regular, then their computational complexity in checking isomorphism can be reduced by considering their disjoint maximal induced regular subgraphs which partition their vertex sets and their interconnected subgraphs. On the other hand, computational complexity on checking graph isomorphism increases with order and degree of the graphs when the graphs are regular, except complete graphs. Circulant graphs are highly symmetric regular graphs and checking of isomorphism among circulant graphs is not easy when the order is large. And once if it is possible to generate isomorphic circulant graphs of a given circulant graph, then it is easy to check their isomorphism. 

A question raised on the isomorphism of $C_{16}(1,3,7)$ and $C_{16}(2,3,5)$ by Elspas and Turner \cite{eltu} was answer by Vilfred \cite{v96} by defining Type-2 isomorphism of $C_n(R)$ w.r.t. $m$ $\ni$ $m$ = $\gcd(n, r) > 1$, $r\in R$ and $r,n\in\mathbb{N}$ and studied such graphs for $m$ = 2 in \cite{v13,v20}. But obtaining Type-2 isomorphic circulant graphs is not easy. Using a $C^{++}$ computer program, the authors obtained families of Type-2 isomorphic $C_{n}(R)$ w.r.t. $m$ = 2,3,5,7 \cite{vw1}-\cite{vw3} and also $C_{np^3}(R)$ w.r.t. $m$ = $p$ \cite{v2-10} for odd prime $p$ and $n\in\mathbb{N}$. In \cite{v2-1}, Vilfred extended the definition of Type-2 isomorphism of circulant graphs $C_n(R)$ w.r.t. $m$ by considering $m > 1$ divides $\gcd(n, r)$, $m^3$ divides $n$ and $r\in R$. Here, $m$ = $\gcd(n, r) > 1$ and $r\in R$ is a particular case of $m > 1$ is a divisor of $\gcd(n, r)$ and $r\in R$. 

In this paper, we present the $C^{++}$ program and also a VB program POLY415.EXE which is used to show how Type-1 and Type-2 isomorphisms take place on a circulant graph as well as for checking and finding Type-1 and Type-2 isomorphic circulant graphs of a given $C_n(R)$. It is also very useful to develop its theory on Type-2 isomorphic circulant graphs \cite{v2-1}-\cite{v2-10}. There are four chapters in this paper. Section 1 is an introductory one and contains basic definitions and a few results which are used in this paper; Section 2 presents the $C^{++}$ computer program that is used to generate Type-2 isomorphic circulant graphs; Section 3 presents the Visual Basic program POLY415.EXE which is used to show and obtain Type-1 and Type-2 isomorphic circulant graphs of a given circulant graph; and Section presents how to use the VB program POLY415.EXE to show and obtain Type-1 and Type-2 isomorphic circulant graphs of a given circulant graph. In \cite{v25} Vilfred extends the definition of Type-2 isomorphism of circulant graphs $C_n(R)$ w.r.t. $m$ by considering $m > 1$ is a divisor of $\gcd(n, r)$ and $r\in R$. In \cite{v2-1}, he further extended the definition by considering $m > 1$ divides $\gcd(n, r)$, $m^3$ divides $n$ and $r\in R$ and used it to develop the theory in \cite{v2-1}-\cite{v2-10}. For basic definitions and results on isomorphic circulant graphs, refer Part 1 \cite{v2-1}.

Through-out this paper, for a set $R$ = $\{ r_1, r_2, \dots, r_k \},$ $C_n(R)$ denotes circulant graph $C_n(r_1, r_2, . . . , r_k)$ where $1 \leq r_1 < r_2 < . . . <$ $r_k \leq [\frac{n}{2}]$. We consider circulant graphs of finite order, $V(C_n(R))$ = $\{v_0, v_1, v_2, . . . , v_{n-1}\}$ with $v_i$ adjacent to $v_{i+r}$ for each $r \in R$, subscript addition taken modulo $n$, $K_n$ = $C_n(1,2,\dots,n-1)$, all cycles have length at least $3$, unless otherwise specified, and the vertices of circulant graph $C_n(R)$ are the corners of a regular $n-gon$, labeled clockwise, $0 \leq i \leq n-1$. However when $\frac{n}{2} \in R$, edge $v_iv_{i+\frac{n}{2}}$ is taken as a single edge for considering the degree of the vertex $v_i$ or $v_{i+\frac{n}{2}}$ and as a double edge while counting the number of edges or cycles in $C_n(R),$ $0 \leq i \leq n-1$. Circulant graphs $C_{16}(1,2,4,7)$, $C_{16}(2,3,4,5)$ and $C_{16}(3,4,5,6)$ are shown in Figures 1, 2 and 3. 

We present here a few definitions and results that are required in this paper.

%%%%%%%%%%%%
% Fig 1,23
%%%%%%%%%%
\begin{figure}[ht]
\centerline{\includegraphics[width=6in]{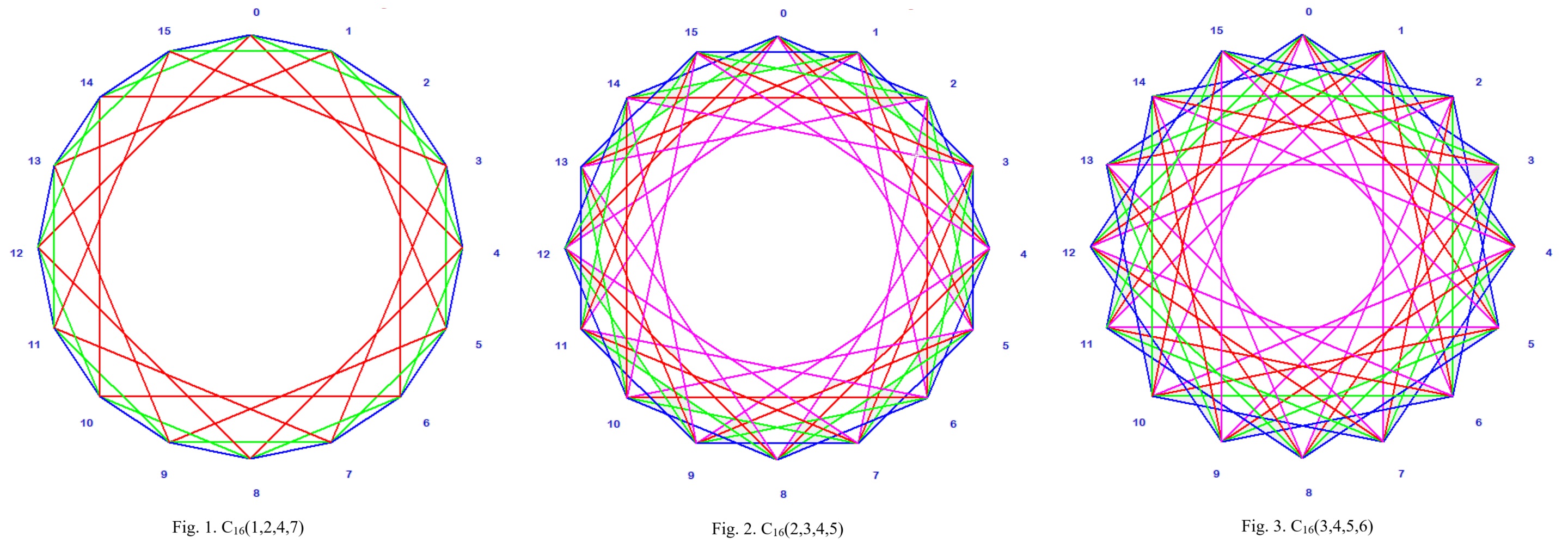}}
\end{figure}
%%%%%%%%%%%%%%%%%%%%%%%%%%%%%%%%%

\begin{definition}{\rm\cite{ad67}} \quad \label{a5} $C_n(R)$ and $C_n(S)$ for $R =$ $\{r_1$, $r_2$, $\ldots$, $r_k\}$ and $S$ = $\{s_1$, $s_2$, $\ldots$, $s_k\}$ are {\em Adam's isomorphic} if there exists a positive integer $x$ relatively prime to $n$ with $S$ = $\{xr_1$, $xr_2$, $\ldots$, $xr_k\}_n^*$ where $<r_i>_n^*$, the {\em reflexive modular reduction} of a sequence $< r_i >$ is the sequence obtained by reducing each $r_i$ modulo $n$ to yield $r_i'$ and then replacing all resulting terms $r_i'$ which are larger than $\frac{n}{2}$ by $n-r_i'$.  

We also call Adam's isomorphism as {\em Type-1 isomorphism}.
\end{definition}

\begin{definition}{\rm \cite{v17}}\quad \label{a6} Let $Ad_n = \{\varphi_ {n,x}: x\in \varphi_ n\},$ $Ad_n(R) = \{\varphi_ {n,x}($ $R): x\in \varphi_ n\} = \{xR:$ $x\in \varphi_ n\}$ and $Ad_n(C_n(R))$ = $T1_n(C_n(R))$ = $\{\varphi_ {n,x}(C_n(R)): x\in \varphi_ n\}$ = $\{C_n(xR): ~x\in \varphi_ n\}$ for a set $R$ = $\{r_1,r_2,...,r_k,n-r_k,n-r_{k-1},...,n-r_1\} \subseteq \mathbb{Z}_n.$ Define $'\circ'$ in $Ad_n(C_n(R))$ such that $\varphi_ {n,x}(C_n(R))$ $\circ$ $\varphi_ {n,y}(C_n(R))$ = $\varphi_ {n,xy}(C_n(R))$ and $C_n(xR)$ $\circ$  $C_n(yR)$ = $C_n((xy)R)$ for every $x,y \in \varphi_ n$, under arithmetic modulo $n$. 

Clearly, $Ad_n(C_n(R))$ is the set of all circulant graphs which are Adam's isomorphic to $C_n(R)$ and $(Ad_n(C_n(R)), \circ )$ = $(T1_n(C_n(R)),~$ $\circ )$ is an Abelian group and we call it as the {\em Adam's group} or {\em Type-1 group of} $C_n(R)$ under $'\circ'$ and we call $Ad_n(C_n(R))$ as the {\em Adam's set} or {\em Type-1 set of} $C_n(R)$.
\end{definition}

$Ad_{16}(C_{16}(1,2,4,7))$ = $\{\varphi_ {16,x}(C_{16}(1,2,4,7)): x\in \varphi_ {16}\}$ = $\{\varphi_ {16,x}(C_{16}(1,2,4,7)): x = 1,3,5,7,9,11,13,15\}$ = $\{C_{16}(1,2,4,7)$, $C_{16}(3,4,5,6)\}$ = $\{\varphi_ {16,x}(C_{16}(1,2,4,7)): x = 1,3\}$. Circulant graphs $C_{16}(1,2,4,7)$, $C_{16}(3,4,5,6)$ which are Adam's isomorphic are given in Figures 1 and 3.  

\begin{definition}{\rm \cite{v2-1}} \quad  \label{d3.4} Let $V(K_n) = \{u_0,u_1,u_2,...,u_{n-1}\}$, $V(C_n(R)) = \{v_0,v_1,v_2,...,$ $v_{n-1}\},$ $r\in R$, $m > 1$ be a divisor of $\gcd(n, r)$ and $|R| \geq 3$.  Define one-to-one mapping $\theta_{n,m,t} :$ $V(C_n(R)) \rightarrow V(K_n)$ such that $\theta_{n,m,t}(v_x)$ = $u_{x+jtm}$,  $\theta_{n,m,t}((v_x, v_{x+s}))$ = $(\theta_{n,m,t}(v_x), \theta_{n,m,t}(v_{x+s}))$ under subscript arithmetic modulo $n$ and $\theta_{n,m,t}(C_n(R))$ = $C_n(\theta_{n,m,t}(R))$ for every $x \in \mathbb{Z}_n$, $s\in R$, $x$ = $qm+j$, $0 \leq j \leq m-1$, $0 \leq q,t \leq \frac{n}{m} -1$ and $\theta_{n,m,t}(R)$ in $C_n(\theta_{n,m,t}(R))$ is calculated under the reflexive modulo $n$. And for a particular value of $t,$ if  $\theta_{n,m,t}(C_n(R))$ = $C_n(S)$ for some $S$  and  $S \neq yR$ for all $y\in \varphi_n$ under reflexive modulo $n$, then $C_n(R)$ and $C_n(S)$ are called {\em isomorphic circulant graphs of Type-2 w.r.t. $m$} and the isomorphism as {\em Type-2 isomorphism w.r.t. $m$.} 
	
When $C_n(R)$ and $C_n(S)$ are Type-2 isomorphic w.r.t. $m$, then we also say that $C_{kn}(kR)$ $(= k.C_n(R))$ and $C_{kn}(kS)$ $(= k.C_n(S))$ are Type-2 isomorphic w.r.t. $m$. 	 
\end{definition}

Let $s \in \mathbb{Z}_n$,  $V_{n,m}$ = $\{\Theta_{n,m,t}:$ $t = 0,1,...,\frac{n}{m}-1\}$, $V_{n,m}(s)$ = $\{\Theta_{n,m,t}(s): t = 0,1,...,\frac{n}{m}-1\}$ and $V_{n,m}(C_n(R))$ = $\{\Theta_{n,m,t}(C_n(R)): t = 0,1,...,\frac{n}{m}-1\}$. Define $'\circ'$ in $V_{n,m}$ such that $\Theta_{n,m,t} ~\circ ~ \Theta_{n,m,t'}$ =  $\Theta_{n,m,t+t'}(\Theta_{n,m,t} ~\circ ~ \Theta_{n,m,t'})(x)$  $( = \Theta_{n,m,t}(\Theta_{n,m,t'}(x))$ = $\Theta_{n,m,t}(x+jt'm)$ = $(x+jt'm)+jtm$ = $x+j(t+t')m )$ = $\Theta_{n,m,t+t'}(x)$ and $\Theta_{n,m,t}(C_n(R)) ~\circ ~ \Theta_{n,m,t'}$ $(C_n(R))$ = $\Theta_{n,m,t+t'}(C_n(R))$ for every $\Theta_{n,m,t},\Theta_{n,m,t'}\in V_{n,m}$ where $t+t'$ is calculated under addition modulo ~$\frac{n}{m}$. Clearly, $(V_{n,m}(s),~ \circ)$ and $(V_{n,m}(C_n(R)),~ \circ)$ are Abelian  groups $\forall s \in \mathbb{Z}_n$. 

\begin{definition}\quad Let $T2_{n,m}(C_n(R))$ = $\{C_n(R)\}$ $\cup$ $\{C_n(S):$ $C_n(S)$ is Type-2 isomorphic of $C_n(R)$ w.r.t. $m\}$ where $r\in R$ and $m > 1$ is a divisor of $\gcd(n, r)$. We call $T2_{n,m}(C_n(R))$ as {\em the Type-2 set of $C_n(R)$
w.r.t. $m$}.

That is, {\em the Type-2 set of $C_n(R)$ w.r.t. $m$} denoted by $T2_{n,m}(C_n(R))$ is $\{C_n(R)\}$ $\cup$ $\{\theta_{n,m,t}(C_n(R)):$ $\theta_{n,m,t}(C_n(R))$ = $C_n(S)$ and $C_n(S)$ is Type-2 isomorphic of $C_n(R)$ w.r.t. $m$, $1 \leq t \leq \frac{n}{m}-1\}$ = $\{\theta_{n,m,0}(C_n(R))\}$ $\cup$ $\{\theta_{n,m,t}(C_n(R)):$ $\theta_{n,m,t}(C_n(R))$ = $C_n(S)$ and $C_n(S)$ is Type-2 isomorphic of $C_n(R)$ w.r.t. $m$, $1 \leq t \leq \frac{n}{m}-1\}$ where $r\in R$ and $m > 1$ is a divisor of $\gcd(n, r)$.
\end{definition}

Clearly, $T2_{n,m}(C_n(R)) \subseteq V_{n,m}(C_n(R))$ and $T1_n(C_n(R))$ $\cap$ $T2_{n,m}(C_n(R))$ = $\{C_n(R)\}$. 

\begin{theorem}{\rm \cite{vc13}}\quad \label{a11} {\rm Let $n \geq 2$, $1 \leq 2s-1 \leq 2n-1$, $n \neq 2s-1$, $R$ = $\{ 2$, $2s-1$, $4n-(2s-1)\}$, $S$ = $\{2$, $2n-(2s-1)$, $2n+2s-1\}$. Then, (i) $\theta_{8n,2,n}(C_{8n}(R))$ = $C_{8n}(S)$ = $\theta_{8n,2,3n}(C_{8n}(R))$, $\theta_{8n,2,n}(C_{8n}(S))$ = $C_{8n}(R)$ = $\theta_{8n,2,3n}(C_{8n}(S))$ and $C_{8n}(R)$ and $C_{8n}(S)$ are Type-2 isomorphic w.r.t. $m$ = 2. \hfill $\Box$}
\end{theorem}

\begin{theorem}{\rm \cite{v2-1}}\quad \label{a12} {\rm Let $n \geq 2$, $k \geq 3$, $1 \leq 2s-1 \leq 2n-1$, $n \neq 2s-1$, $R$ = $\{ 2s-1$, $4n-(2s-1)$, $2p_1$, $2p_2$, $\dots$, $2p_{k-2} \}$, $S$ = $\{2n-(2s-1)$, $2n+2s-1$, $2p_1,2p_2,\dots,2p_{k-2}\}$, $2y\in R,S$, $\gcd(4n,y)$ = 1, $p_1,p_2,\dots,p_{k-2} \in \mathbb{N}$ and $\gcd(p_1,p_2,\dots,p_{k-2})$ = 1. Then, (i) $\theta_{8n,2,n}(C_{8n}(R))$ = $C_{8n}(S)$ = $\theta_{8n,2,3n}(C_{8n}(R))$, $\theta_{8n,2,n}(C_{8n}(S))$ = $C_{8n}(R)$ = $\theta_{8n,2,3n}(C_{8n}(S))$ and (ii) for a given set of values of $p_1,p_2,\dots,p_{k-2}$ and $y$, $C_{8n}(R)$ and $C_{8n}(S)$ are isomorphic of either Adam's or Type-2 w.r.t. $m$ = 2. Moreover, for all such possible values of $p_1,p_2,\dots,p_{k-2}$ and $y$, the set $\{ C_n(S) = \theta_{8n,2,n}(C_{8n}(R)): p_1,p_2,\dots,p_{k-2} \in \mathbb{N}\}$ contains all Type-2 isomorphic circulant graphs of $C_{8n}(R)$ w.r.t. $m$ = 2.  \hfill $\Box$}
\end{theorem}

\begin{rem}{\rm \cite{v2-1}}\quad \label{a12} Under the transformation $\Theta_{n,m,t}$ acting on $C_n(R)$, values which are integer multiple of $m$ are not changing where $m > 1$ is a divisor of $gcd(n, r)$ and $r\in\mathbb{N}$. Also, for a given $C_n(R)$, with respect to different $m$ or $r$ or both, we may get different Type-2 isomorphic circulant graphs.
\end{rem}

\begin{theorem}  {\rm \cite{v2-3}}\label{b1} Let $C_n(R)$ be a circulant graph, $m > 1$ be a divisor of $gcd(n, r)$ and $r\in R$. Let $t_1$ be the smallest positive integer $\ni$ $\Theta_{n,m,t_1}(C_n(R))$ = $C_n(S)$ for some $S$ and $C_n(R)$ and $C_n(S)$ be Type-2 isomorphic w.r.t. $m$. Then, $T2_{n,m}(C_n(R))$ = $\{\Theta_{n,m,jt_1}(C_n(R)):~ j = 0,1,2,...,q-1$, $mqt_1 = n\}$ and $(T2_{n,m}(C_n(R)), \circ)$ is an Abelian subgroup of $(V_{n,m}(C_n(R)), \circ)$ where $\Theta_{n,m,t}(C_n(R)) \circ \Theta_{n,m,t'}(C_n(R))$ = $\Theta_{n,m,t+t'}(C_n(R))$, $0 \leq t,t' \leq$ $\frac{n}{m}-1$ and $t+t'$ is calculated under arithmetic modulo $\frac{n}{m}$. \hfill $\Box$
\end{theorem}

\begin{definition}{\rm \cite{v2-3}}\quad \label{b2} With usual notation, group $(T2_{n,m}(C_n(R)), \circ)$ is called {\em the Type-2 group of $C_n(R)$ w.r.t. $m$} under $~\lq\circ\rq.$  
\end{definition}

\begin{theorem} {\rm \cite{v2-6}} \label{d9} Let $p$ be an odd prime number, $1 \leq x \leq p-1$,  $1 \leq x+yp \leq np^2-1$, $1 \leq i,j \leq p$, $y\in\mathbb{N}_0$, $d^{np^3, x+yp}_i$ = $(i-1)npx+$ $x+yp$, $R^{np^3, x+yp}_i$ = $\{p$, $d^{np^3, x+yp}_i$, $np^2-d^{np^3, x+yp}_i$, $np^2+d^{np^3, x+yp}_i$, $2np^2-$ $d^{np^3, x+yp}_i$, $2np^2+d^{np^3, x+yp}_i$, $3np^2-d^{np^3, x+yp}_i$, $3np^2+d^{np^3, x+yp}_i$, . . . , $(p-1)np^2-d^{np^3, x+yp}_i$, $(p-1)np^2+d^{np^3, x+yp}_i$, $np^3-d^{np^3, x+yp}_i$, $np^3-p\}$ and $i,j,n,x\in \mathbb{N}$.  Then, for given $n, x, y$ and $p$, $\Theta_{np^3,p,jn} (R^{np^3, x+yp}_i)$ =  $R^{np^3, x+yp}_{i+j}$, $\Theta_{np^3,p,jn}(C_{np^3}(R^{np^3, x+yp}_i))$ = $C_{np^3}(R^{np^3, x+yp}_{i+j})$ and the $p$ circulant graphs $C_{np^3}(R^{np^3, x+yp}_i)$ for $i = 1,2,...,p$ are Type-2 isomorphic w.r.t. $m$ = $p$ where $i+j$ in $R^{np^3, x+yp}_{i+j}$ is calculated under addition modulo $p$. \hfill $\Box$
\end{theorem}

\section{Our $C^{++}$ computer program to generate Type-2 isomorphic circulant graphs}

The following $C^{++}$ computer program \cite{bs} is used to generate Type-2 isomorphic circulant graphs and from these we obtained Theorems \ref{a11}, \ref{a12}, \ref{b1} and \ref{d9}.

\vspace{.5cm}
\noindent
{\bf Our $C^{++}$ computer program}

\vspace{.2cm}
\noindent
{\bf Program description}

This program finds Type-2 isomorphic circulant graphs of any given circulant graph of finite order. It reads a prime number $‘p’$ as the input. Here, we consider different values of $‘n’$ from 1 to 5, then the number of vertices of the graph $‘v’$ is computed as $v$ = $np^3$. For a given $‘p’$ and $‘v’$, different possible jump sizes of circulant graphs are considered and computed to find Type-2 isomorphic circulant graphs, if exist. The function ‘ComputeTrans’, finds $‘p’$ sets of jump sizes corresponding to the $‘p’$ transformations $\Theta_{np^3,p,jn}(C_{np^3}(R^{np^3, x+yp}_i))$, $j$ = 1 to $p$ for the given circulant graph which are Type-2 isomorphic to the given graph, see Theorem \ref{d9}. The program is written in $C^{++}$ as follows.

\vspace{.2cm}
\noindent
{\bf Program}

\vspace{.2cm}
\noindent
$\#include<iostream.h>$
\\
$\#include<iomanip.h>$
\\
$\#include<conio.h>$

\vspace{.2cm}
\noindent
void main()
\\
$\{$
\\
 $void ~ComputeTrans(int^*, int, int, int, int);$
 \\
void DispNums($int^*$,int,int);
\\
 void Sort($int^*$, int);
\\
 $int~ n,p,v,jumpsize,js,i,j,k$;
\\
 $clrscr()$;

 $cout<<``Enter~ p..."$;

 $cin>>p$;
\\
 for($n$=1;$n<=5$; $n++$)
\\
 $\{$

 	$cout<<``$ $\ n$ $p$ = $"<< p <<``$ and $n="<< n <<``\ n"$;

 	$cout<<``----------------------------\ n"$;

 	$v$ = $n\times p \times p \times p$;

 	jumpsize = $p+1$;

 	js = $2\times jumpsize$;

 	$int^*ptr$ = $new~ int[js]$;

  	for($k$=1; $k<=n\times p-1$;$k++$)

 	$\{$

  \hspace{.5cm}     if $(k ~\% ~p == 0)$

  \hspace{.5cm}     	continue;

 \hspace{.75cm}		$j$ =1;

 \hspace{.75cm}		$ptr[0]$ = $k$;

 \hspace{.75cm}		$ptr[1]$ = $p$;

 \hspace{.75cm}		for($i$ =2; $i < jumpsize$; $i+$ = 2)

 	 \hspace{.5cm}	$\{$
  			
  \hspace{.75cm}	         $ptr[i]$ = $j\times n\times p \times p-k$;

  \hspace{.75cm}		$ptr[i+1]$ = $j\times n\times p\times p+k$;

  \hspace{.75cm}		$j++$;
 		
           \hspace{.5cm}    $\}$
 		
    \hspace{.5cm}                for($i$ =0; $i < jumpsize$; $i++$)

  \hspace{.5cm}		$\{$
  			
  \hspace{.75cm}	         $ptr[2\times jumpsize-1-i]$ = $v-ptr[i]$;
 		 
   \hspace{.5cm}                  $\}$

   \hspace{.5cm}    	$Sort(ptr,js)$;

  \hspace{.5cm}	$DispNums(ptr,v,js)$;

  \hspace{.5cm}	$ComputeTrans(ptr,v,p,n,js)$;

  \hspace{.5cm}	$\}$
 
 \hspace{.5cm}     $\}$
 
 \hspace{.5cm}    $getch()$;

$\}$

\vspace{.2cm}
\noindent
void $DispNums(int^*ptr,int v,int js)$

$\{$
	
\hspace{.5cm}	 $cout<<"C"<<v<<"("$;

\hspace{.5cm}	 for(int $i$ =0; $i<js-1$; $i++$)

\hspace{.5cm}	$cout<<ptr[i]<<","$;

\hspace{.5cm}	 $cout<<ptr[js-1]<<")\ n"$;

$\}$

\vspace{.2cm}
\noindent
void $ComputeTrans(int^*ptr, int v, int p, int n, int js)$

$\{$

	void $Sort(int *, int)$;

	void $DispNums(int *, int, int)$;

	$int i,t,jmp,nt$;

	$int^*$ NewNums;

	NewNums=new $int[js]$;

	$nt=p$;

   for($t$ =1; $t<nt$; $t++$)

	$\{$

		for($i$ =0; $i<js$; $i++$)

\hspace{.5cm}	$\{$

\hspace{.75cm}		$jmp=ptr[i]\% p$;

\hspace{.75cm}		$NewNums[i]=(ptr[i] + n*t*p*jmp)\% v$;

 \hspace{.5cm}	 	$\}$

		$Sort(NewNums,js)$;

		$DispNums(NewNums,v,js)$;

	$\}$
  
$cout<<``\ n----------------------------------\ n\ n"$;

$\}$

  // Sorts an One-dimesional array

\vspace{.2cm}
\noindent
void $Sort(int^*ptr, int js)$

  $\{$

 \hspace{.5cm}	$int i,j,tmp$;

 \hspace{.5cm}	for($i$ =0; $i<js-1$; $i++$)

 \hspace{.5cm}	$\{$

\hspace{.75cm}		 for($j$ = $i+1$; $j<js$; $j++$)

\hspace{.75cm}		$\{$

\hspace{1cm}		  if ($ptr[i]>ptr[j]$)

\hspace{1.5cm}		  $\{$

\hspace{1.5cm}			$tmp=ptr[i]$;

\hspace{1.5cm}			$ptr[i]=ptr[j]$;

\hspace{1.5cm}			$ptr[j]=tmp$;

\hspace{1cm}		  $\}$
	 
\hspace{.75cm}	           $\}$

 \hspace{.5cm}	$\}$

  $\}$

\vspace{.3cm}
In \cite{v24}, a list of Type-2 sets $T2_{np^3,p}(C_{np^3}(R^{np^3,x+yp}_i))$ of isomorphic circulant graphs of Type-2 w.r.t. $p$, generated by the above program, is given in the Annexure for $p$ = 3,5,7,11 and $n$ = 1 to 5 and also for $p$ = 13 and $n$ = 1 to 3 where $(T2_{np^3,p}(C_{np^3}(R^{np^3,x+yp}_i)), \circ)$ is a subgroup of $(V_{np^3,p}(C_{np^3}(R^{np^3,x+yp}_i)), \circ)$, $1 \leq i \leq p$, $1 \leq x \leq p-1$, $1 \leq x + yp \leq np^2-1$, $0 \leq y \leq np-1$ and $p,np^3-p\in R^{np^3,x+yp}_i$. 

\section{Visual Basic Computer Program POLY415.EXE} 

Visual Basic (VB) programs \cite{hv} are well-suited for creating user-friendly interfaces and interacting with databases. The following VB program POLY415.EXE demonstrates the movements of vertices of circulant graph $C_n(R)$ and displays how Type-1 and Type-2 isomorphisms of a given circulant graph is taking place. It also displays Type-1 and Type-2 isomorphic circulant graphs of a given circulant graph, if they exist. 

\vspace{.5cm}
\noindent
{\bf Visual Basic Program POLY415.EXE}

\vspace{.3cm}
\noindent
Private Sub $Command1\_Click()$
\\
Dim $i$ As Integer, $j$ As Integer
\\
$n$ = $Int(Val(Text1.Text))$
\\
If List1.SelCount $<=$ 0 Then

 MsgBox (``Select Atleast One Jump value")

 List1.SetFocus
\\
End If

JumpSize = List1.SelCount

Cycles = Int(Val(Text2.Text))

ReDim JumpVal(JumpSize)

$j$ = 0
\\
For $i$ = 0 To List1.ListCount - 1
 
If $List1.Selected(i)$ = True Then

 $JumpVal(j)$ = $List1.List(i)$

$j$ = $j + 1$

 End If
\\
Next $i$

Call AdjacencyMatrix(JumpVal)

Call Form2.InitGraph
\\
End Sub

\vspace{.2cm}
\noindent
Private Sub $Text1\_LostFocus()$
\\
Dim $i$ As Integer, num As Integer
\\
$n$ = Int(Val(Text1.Text))
\\
num = $Int(n / 2)$
\\
List1.Clear
\\
For $i$ = 1 To num

List1.AddItem i
\\
Next $i$
\\
End Sub

\vspace{.2cm}
\noindent
Dim $k$ As Integer, $s$ As Integer, NumMoves As Integer, $nm$ As Integer
\\
Dim $X(), Y(), X1(), Y1(), X2(), Y2(), OldX(), OldY(), OldX1(), OldY1(), OldTheta(), NewTheta()$
\\
Dim $MidX, MidY$, Radius As Single

\vspace{.2cm}
\noindent
Sub $InitGraph()$
\\
Dim $i$ As Integer
\\
Dim StartTheta As Double, Theta As Double, StepVal As Double
\\
ReDim $X(n), Y(n), X1(n), Y1(n), X2(n), Y2(n), OldX(n), OldY(n), OldX1(n), OldY1(n), OldTheta(n), NewTheta(n)$
\\
ReDim $Labels(n)$
\\
Form2.Show

Form2.FillStyle = 0

MidX = Form2.ScaleWidth / 2

MidY = Form2.ScaleHeight / 2

Radius = Form2.ScaleHeight / 2.5

StartTheta = $3\times \Pi / 2$

Form2.DrawWidth = 2

StepVal = $2\times \Pi / n$

$\theta$ = StartTheta
\\
For $i$ = 0 To $n-1$

$X(i)$ = $MidX + Radius \times \cos(\theta)$

$Y(i)$ = $MidY + Radius \times \sin(\theta)$

$X1(i)$ = $MidX + (Radius + 520) \times \cos(\theta)$

$Y1(i)$ = $MidY + (Radius + 520) \times \sin(\theta)$

$X2(i)$ = $MidX + (Radius + 150) \times \cos(\theta)$

$Y2(i)$ = $MidY + (Radius + 150) \times \sin(\theta)$

$OldX(i)$ = $X(i)$

$OldY(i)$ = $Y(i)$

$OldX1(i)$ = $X1(i)$

$OldY1(i)$ = $Y1(i)$

$OldTheta(i)$ = $Theta$

$NewTheta(i)$ = $OldTheta(i)$

$\theta$ = $\theta + StepVal$
\\
Next $i$
\\
$Form2.Labels(0).Caption$ = 0

$Form2.Labels(0).Left$ = $X1(0)$

$Form2.Labels(0).Top$ = $Y1(0)$

$Form2.Labels(0).Visible$ = True
\\
For $i$ = 1 To $n-1$

$Load Form2.Labels(i)$
\\
Next $i$

\vspace{.2cm}
Call $DrawCirculant(X, Y)$

Call $PutLabel(X1, Y1)$
\\
End Sub

\vspace{.2cm}
\noindent
Sub $DrawCirculant(X, Y)$
\\
Dim $i$ As Integer, $j$ As Integer, $k$ As Integer
\\
For $i$ = 0 To $n-1$
\\
For $j$ = i To $n-1$
\\
If $AdMat(i, j) <>$ 0 Then

 Form2.Line $(X(i), Y(i))-(X(j), Y(j)), AdMat(i, j)$
\\
End If
\\
Next $j$
\\
Next $i$
\\
End Sub

\vspace{.2cm}
\noindent
Sub $HideCirculant(X, Y)$
\\
Dim $i$ As Integer, $j$ As Integer, $k$ As Integer
\\
For $i$ = 0 To $n-1$
\\
For $j$ = $i$ To $n-1$
\\
If $AdMat(i, j) <>$ 0 Then

 Form2.Line $(X(i), Y(i))-(X(j), Y(j))$, Form2.BackColor
\\
End If
\\
Next $j$
\\
Next $i$
\\
End Sub

\vspace{.2cm}
\noindent
Public Sub $Rotate()$
\\
Dim steps As Single
\\
Dim incr As Double
\\
Dim $i, m$ As Integer

steps = 20

incr = $2 \times \Pi / (n \times steps)$

Call $HideCirculant(X, Y)$
\\
For $m$ = 1 To Cycles - 1
\\
For $i$ = $m$ To $n-1$ Step Cycles

$NewTheta(i)$ = $NewTheta(i) + m \times incr$

$X(i)$ = $MidX + Radius \times \cos(NewTheta(i))$

$Y(i)$ = $MidY + Radius \times \sin(NewTheta(i))$

$X1(i)$ = $MidX + (Radius + 520) \times \cos(NewTheta(i))$

$Y1(i)$ = $MidY + (Radius + 520) \times \sin(NewTheta(i))$
\\
Next $i$
\\
Next $m$

Call $DrawCirculant(X, Y)$

Call $PutLabel(X1, Y1)$

Call Delay
\\
End Sub

\vspace{.2cm}
\noindent
Private Sub $Delay()$
\\
Dim Start
\\
Dim Check

    Start = Timer
\\
 Do Until Check $>= Start + ((Form2.HScroll1.Max - Form2.HScroll1.Value) / (1000 \times Cycles))$

        Check = Timer

    Loop
\\
End Sub

\vspace{.2cm}
\noindent
Private Sub $HideLabel()$
\\
Dim $i$ As Integer
\\
For $i$ = 0 To $n-1$

$Form2.Labels(i).Visible$ = False
\\
Next $i$
\\
End Sub

\vspace{.2cm}
\noindent
Private Sub $PutLabel(X1, Y1)$
\\
Dim $i$ As Integer
\\
For $i$ = 0 To $n-1$

$Form2.Labels(i).Caption$ = i

$Form2.Labels(i).Left$ = $X1(i)$

$Form2.Labels(i).Top$ = $Y1(i)$

$Form2.Labels(i).Visible$ = True
\\
Next $i$
\\
End Sub

\vspace{.2cm}
\noindent
Public Sub AdamRotate($k$ As Integer)
\\
Dim steps As Single
\\
Dim incr As Double

steps = 15

incr = $2 \times \Pi / (n \times steps)$

$j\%$ = 1

Call $HideCirculant(X, Y)$
\\
Do While $j$ $<=$ steps
\\
For $i\%$ = 1 To $n-1$

$NewTheta(i)$ = $NewTheta(i) + i \times (k - 1) \times incr$

$X(i)$ = $MidX + Radius \times \cos(NewTheta(i))$

$Y(i)$ = $MidY + Radius \times \sin(NewTheta(i))$

$X1(i)$ = $MidX + (Radius + 520) \times \cos(NewTheta(i))$

$Y1(i)$ = $MidY + (Radius + 520) \times \sin(NewTheta(i))$
\\
Next $i$

Call $DrawCirculant(X, Y)$

Call $PutLabel(X1, Y1)$

Call Delay

Call $HideCirculant(X, Y)$

$j$ = $j + 1$

Loop

Call $DrawCirculant(X, Y)$
\\
End Sub

\vspace{.2cm}
\noindent
Sub $Reset()$

Call $HideCirculant(X, Y)$

Call $DrawCirculant(OldX, OldY)$

Call $PutLabel(OldX1, OldY1)$
\\
For $i\%$ = 0 To $n-1$

$X(i)$ = $OldX(i)$

$Y(i)$ = $OldY(i)$

$X1(i)$ = $OldX1(i)$

$Y1(i)$ = $OldY1(i)$

$NewTheta(i)$ = $OldTheta(i)$
\\
Next $i$
\\
End Sub

\vspace{.2cm}
\noindent
Sub $StepMove()$
\\
Dim steps As Single
\\
Dim incr As Double
\\
Dim $i$ As Integer, $m$ As Integer, $nsteps$ As Integer

steps = 20

incr = $2 \times \Pi / (n \times steps)$

Call $HideCirculant(X, Y)$

$nsteps$ = 0
\\
Do While $nsteps < (Cycles \times steps)$
\\
For $m$ = 1 To $Cycles - 1$
\\
For $i$ = $m$ To $n-1$ Step Cycles

$NewTheta(i)$ = $NewTheta(i) + m \times incr$

$X(i)$ = $MidX + Radius  \times \cos(NewTheta(i))$

$Y(i)$ = $MidY + Radius  \times \sin(NewTheta(i))$

$X1(i)$ = $MidX + (Radius + 520)  \times \cos(NewTheta(i))$

$Y1(i)$ = $MidY + (Radius + 520)  \times \sin(NewTheta(i))$

Next $i$

Next $m$

$nsteps$ = $nsteps + 1$

Call $DrawCirculant(X, Y)$

Call $PutLabel(X1, Y1)$

Call Delay

Call $HideCirculant(X, Y)$

Loop

Call $DrawCirculant(X, Y)$
\\
End Sub

\vspace{.2cm}
\noindent
Private Sub $StepMoveButtn\_Click()$
\\
$StepMoveButtn.Visible$ = False

$AdamIsoButtn.Visible$ = False

$RotateButtn.Visible$ = False

$StopButtn.Visible$ = False

$ContinueButtn.Visible$ = False

$NextMoveButtn.Visible$ = True

$ResetButtn.Visible$ = True

$Label2.Visible$ = False

$Label5.Visible$ = False

$AdLabel.Visible$ = True

$Label7.Caption$ = ``Transformation"

Call Reset

Call InitJumpVal

Call RelPrime

Call ComputeAdamIso

$s$ = 1

Call StepMove

Call CompJumpVal

$Label6.Caption$ = $s \&$ ``."

NumMoves = $n / Cycles$
\\
End Sub

\vspace{.2cm}
\noindent
Private Sub $AdamIsoButtn\_Click()$

HScroll1.Visible = False

HScroll1.Enabled = False

Label2.Visible = False

Label5.Visible = True

AdLabel.Visible = False

StepMoveButtn.Visible = False

AdamIsoButtn.Visible = False

RotateButtn.Visible = False

StopButtn.Visible = False

NextMoveButtn.Visible = False

ContinueButtn.Visible = True

Label7.Caption = ``ADAM Isomorphic Graphs"

Call RelPrime

$k$ = 1

$nm$ = 1

Call Reset

Call $AdamRotate(RP(k))$

Label6.Caption = $nm \&$ ``."

Label5.Caption = ``$m$ =" $\& RP(k)$

Call $ComputeJV(RP(k))$
\\
End Sub

\vspace{.2cm}
\noindent
Private Sub $RotateButtn\_Click()$

HScroll1.Visible = True

HScroll1.Enabled = True

Label2.Visible = True

Label5.Visible = False

AdLabel.Visible = False

Label3.Visible = False

Label6.Visible = False

Label7.Caption = `` "

StepMoveButtn.Visible = False

AdamIsoButtn.Visible = False

RotateButtn.Visible = True

StopButtn.Visible = True

NextMoveButtn.Visible = False

ContinueButtn.Visible = False

ResetButtn.Visible = True

Timer1.Enabled = True
\\
End Sub

\vspace{.2cm}
\noindent
Private Sub $StopButtn\_Click()$

HScroll1.Enabled = False

HScroll1.Visible = False

Label2.Visible = False

Timer1.Enabled = False
\\
End Sub

\vspace{.2cm}
\noindent
Private Sub $ContinueButtn\_Click()$

$k$ = $k + 1$
\\
If $k >= RPCnt$ Then

ContinueButtn.Visible = False

StepMoveButtn.Visible = True

AdamIsoButtn.Visible = True

RotateButtn.Visible = True

StopButtn.Visible = False

NextMoveButtn.Visible = False
\\
Else

Call Reset

Call $AdamRotate(RP(k))$

$nm$ = $nm + 1$

Label6.Caption = $nm \&$ ``."

Label5.Caption = ``$m$ =" $\& RP(k)$

Call $ComputeJV(RP(k))$
\\
End If
\\
End Sub

\vspace{.2cm}
\noindent
Private Sub $NextMoveButtn\_Click()$

$s$ = $s + 1$

Call StepMove

Call CompJumpVal

Label6.Caption = $s \&$ ``."
\\
If $s > NumMoves - 1$ Then

NextMoveButtn.Visible = False

StepMoveButtn.Visible = True

AdamIsoButtn.Visible = True

RotateButtn.Visible = True

StopButtn.Visible = False

ContinueButtn.Visible = False
\\
End If
\\
End Sub

\vspace{.2cm}
\noindent
Private Sub $ResetButtn\_Click()$

StepMoveButtn.Visible = True

AdamIsoButtn.Visible = True

StopButtn.Visible = False

ResetButtn.Visible = False

Call Reset
\\
End Sub

\vspace{.2cm}
\noindent
Private Sub $ExitButtn\_Click()$

Form1.Hide

Form2.Hide
\\
End
\\
End Sub

\vspace{.2cm}
\noindent
Private Sub $Form\_Load()$

Form2.DrawMode = 13

Timer1.Enabled = False

StopButtn.Visible = False

ContinueButtn.Visible = False

NextMoveButtn.Visible = False

ResetButtn.Visible = False

Label2.Visible = False

Call $DispTitle(Label4, n, JumpSize, JumpVal)$
\\
End Sub

\vspace{.2cm}
\noindent
Private Sub $Timer1\_Timer()$

Call Rotate
\\
End Sub

\vspace{.2cm}
\noindent
Public $AdMat()$
\\
Public Const $\Pi$ As Double = 3.1415926535
\\
Public $n$ As Integer
\\
Public Cycles As Integer
\\
Public JumpSize As Integer
\\
Public RPCnt As Integer
\\
Public $JumpVal(), NewJV(), NewJumpVal()$
\\
Public $RP()$ As Integer
\\
Public $AdamMat()$

Public Function GCD($a$ As Integer, $b$ As Integer) As Integer
\\
If $a$ Mod $b$ = 0 Then

GCD = $b$
\\
Else

GCD = GCD($b$, $a$ Mod $b$)
\\
End If
\\
End Function

\vspace{.2cm}
\noindent
Public Sub DispTitle(Lbl As Control, $m$ As Integer, $JS$ As Integer, $JV$)
\\
Dim Title As String

Title = ``$C$" $\&$ $m$ $\&$ ``("

For $i\%$ = 0 To $JS - 2$

Title = Title $\&$ $JV(i)$ $\&$ ``,"

Next $i$

Title = Title $\&$ $JV(JS - 1)$ $\&$ ``)"

Lbl.Caption = Title
\\
End Sub

\vspace{.2cm}
\noindent
Public Sub ComputeJV($m$ As Integer)
\\
ReDim NewJV(JumpSize)
\\
For $i\%$ = 0 To JumpSize - 1

$NewJV(i)$ = $JumpVal(i) \times m$

$NewJV(i)$ = $NewJV(i) Mod ~n$
\\
If $NewJV(i) > Int(n / 2)$ Then
 
$NewJV(i)$ = $n - NewJV(i)$
\\
 End If
\\
Next $i$

Call Sort(JumpSize, NewJV)

Call DispTitle(Form2.Label3, $n$, JumpSize, NewJV)
\\
End Sub

\vspace{.2cm}
\noindent
Public Sub Sort(m As Integer, JV)
\\
For $i\%$ = 0 To $m - 2$
\\
 For $j\%$ = $i$ To $m - 1$
\\
  If $JV(i) > JV(j)$ Then

   temp = $JV(i)$

   $JV(i)$ = $JV(j)$

   $JV(j)$ = temp
\\
  End If
\\
 Next $j$
\\
 Next $i$
\\
End Sub

\vspace{.2cm}
\noindent
Public Sub $CompJumpVal()$
\\
Dim DJumpSize As Integer
\\
Dim IsSym As Boolean
\\
Dim IsAdam As Boolean
\\
DJumpSize = $2 \times JumpSize$
\\
For $i\%$ = 0 To $DJumpSize - 1$

$NewJumpVal(i)$ = $NewJumpVal(i) + (NewJumpVal(i) Mod Cycles) \times Cycles$
\\
If $NewJumpVal(i) > n$ Then

 $NewJumpVal(i)$ = $NewJumpVal(i)~ Mod ~n$
\\
End If
\\
Next $i$

Call Sort(DJumpSize, NewJumpVal)

IsSym = CheckSym(DJumpSize, NewJumpVal)
\\
If IsSym = True Then

Call DispTitle(Form2.Label3, $n$, JumpSize, NewJumpVal)

IsAdam = CheckAdamIso(JumpSize, NewJumpVal)
\\
If IsAdam = True Then

Form2.AdLabel = ``Adams"
\\
Else

Form2.AdLabel.Caption = ``Non-Adams"
\\
End If
\\
Else

 Form2.Label3.Caption = ``Non-Circulant"

 Form2.AdLabel.Caption = `` "
\\
End If
\\
End Sub

\vspace{.2cm}
\noindent
Public Sub $InitJumpVal()$
\\
Dim DJumpSize As Integer
\\
DJumpSize = $2 \times JumpSize$
\\
ReDim NewJumpVal(DJumpSize)
\\
For $i\%$ = 0 To $JumpSize - 1$

$NewJumpVal(i)$ = $JumpVal(i)$

$NewJumpVal(DJumpSize - (i + 1))$ = $n - JumpVal(i)$
\\
Next $i$
\\
End Sub

\vspace{.2cm}
\noindent
Function CheckSym($m$ As Integer, JV) As Boolean

CheckSym = True
\\
For $i\%$ = 0 To $JumpSize - 1$
\\
If $JV(i) + JV(m - (i + 1)) <> n$ Then

 CheckSym = False

  Exit For
\\
End If
\\
Next $i$
\\
End Function

\vspace{.2cm}
\noindent
Public Sub $RelPrime()$
\\
Dim $m$ As Integer

ReDim $RP(n)$

$RPCnt$ = 0
\\
For $i\%$ = 1 To $n - 1$

$m$ = $\gcd(n, i)$
\\
If $m$ = 1 Then

 $RP(RPCnt)$ = $i$

 $RPCnt$ = $RPCnt + 1$
\\
End If
\\
Next $i$
\\
End Sub

\vspace{.2cm}
\noindent
Public Sub $ComputeAdamIso()$

ReDim AdamMat(RPCnt, JumpSize)

ReDim NewJV(JumpSize)
\\
For $i\%$ = 0 To $RPCnt / 2$
\\
For $j\%$ = 0 To $JumpSize - 1$

   $NewJV(j)$ = $JumpVal(j) \times RP(i)$

   $NewJV(j)$ = $NewJV(j)~ Mod ~n$
\\
   If $NewJV(j) > Int(n / 2)$ Then

     $NewJV(j)$ = $n - NewJV(j)$
\\
   End If
\\
 Next $j$

 Call Sort(JumpSize, NewJV)
\\
 For $j$ = 0 To $JumpSize - 1$

    $AdamMat(i, j)$ = $NewJV(j)$
\\
 Next $j$
\\
Next $i$
\\
End Sub

\vspace{.2cm}
\noindent
Function CheckAdamIso($m$ As Integer, JV) As Boolean
\\
Dim Status As Boolean
\\
For $i\%$ = 0 To $RPCnt - 1$

Status = True
\\
For $j\%$ = 0 To $m - 1$
\\
If $AdamMat(i, j) <> JV(j)$ Then

Status = False

Exit For
\\
End If
\\
Next $j$

If Status = True Then

Exit For
\\
End If
\\
Next $i$

CheckAdamIso = Status
\\
End Function

\vspace{.2cm}
\noindent
Sub AdjacencyMatrix(JumpVal)
\\
Dim $i, j$
\\
Dim color, clr

ReDim $AdMat(n, n)$
\\
For $i$ = 0 To $n - 1$

    $AdMat(0, i)$ = 0
\\
Next $i$

clr = 9
\\
For $j$ = 0 To $JumpSize - 1$

        color = $QBColor(clr)$

        $AdMat(0, JumpVal(j))$ = color

        $AdMat(0, n - JumpVal(j))$ = color

        $clr$ = $clr + 1$

        If $clr$ = 11 Then

         $clr$ = $clr + 1$

        End If
\\
Next j
\\
For $i$ = 1 To $n - 1$
\\
    For $j$ = 0 To $n - 1$

        $AdMat(i, j)$ = $AdMat(0, (n + j - i) Mod n)$
\\
    Next $j$
\\
Next $i$
\\
End Sub

\section{How to use the program POLY415.EXE}

The Visual Basic program execution file POLY415.EXE is used to show how Type-1 and Type-2 isomorphisms of a given circulant graph is taking place. Here, we present the steps to operate POLY415.EXE and  demonstrate with an example. We consider the following list of Type-2 isomorphic circulant graphs and choose any one for demonstration. List of more Type-2 isomorphic circulant graphs are given in \cite{v2-6}. One is free to choose any circulant graph.  

\begin{enumerate} 
\item [\rm (i)] For $p$ = 3, $x$ = 1, $y$ = 0 and $n$ = 1,
\\
$T2_{np^3,p}(C_{np^3}(R^{np^3,x+yp}_i))$ 

\hspace{1cm} = $T2_{27,3}(C_{27}(R^{27,1}_i))$ = $\{C_{27}(1,3,8,10)$, $C_{27}(3,4,5,13)$, $C_{27}(2,3,7,11)\}$;

\item [\rm (ii)]  For $p$ = 3, $x$ = 2, $y$ = 1 and $n$ = 2, 
\\
$T2_{np^3,p}(C_{np^3}(R^{np^3,x+yp}_i))$ 

\hspace{1cm} = $T2_{54,3}(C_{54}(R^{54,5}_i))$  = $\{C_{54}(3,5,13,23)$, $C_{54}(1,3,17,19)$, $C_{54}(3,7,11,25)\}$;

\item [\rm (iii)] For $p$ = 3, $x$ = 2, $y$ = 2 and $n$ = 3, 
\\
$T2_{np^3,p}(C_{np^3}(R^{np^3,x+yp}_i))$ 

\hspace{1cm} = $T2_{81,3}(C_{81}(R^{81,8}_i))$ = $\{C_{81}(3,8,19,35)$, $C_{81}(1,3,26,28)$, $C_{81}(3,10,17,37)\}$.
\end{enumerate}

\vspace{.2cm}
\noindent
{\bf Steps used to activate the program POLY415.EXE}

\vspace{.2cm}
The following steps are used to activate the program POLY415.EXE. Suppose we want to use the program for the circulant graph $C_{27}(1,3,8,10)$ which is chosen from the list given above.

\begin{enumerate} 
\item [\rm (i)] Double click POLY415.EXE to open the file. In the screen, we get as in Figure 4.

\item [\rm (ii)] Enter the data in the displayed window as 

27 in Vertices; Select 1, 3, 8, 10 as jump size values and 3 (= $m$) in cycles. Then the window looks as in Figure 5.

\item [\rm (iii)] Now, click OK button. The circulant graph $C_{27}(1,3,8,10)$ is displayed in the screen as shown in Figure 6.

\item [\rm (iv)] Choose StepMove or Rotate to see how the graph changes under the transformation used to define Type-2 isomorphism. Repeat this step to get Type-2 isomorphic circulant graphs, if exist. 

When StepMove button is clicked, $t$ takes the value 1 and circulant graph $C_{27}(3,4,5,13)$ corresponding to $\theta_{27,3,t=1}(C_{27}(1,3,8,10))$ is displayed as shown in Figure 7. Moreover, it indicates  $C_{27}(3,4,5,13)$ = $\theta_{27,3,1}(C_{27}(1,3,8,10))$ is Type-2 isomorphic w.r.t. $m$ = 3 (to $C_{27}(1,3,8,10)$).

When we click the NextMove button, $t$ takes the value 2 and circulant graph $C_{27}(2,3,7,11)$ corresponding to $\theta_{27,3,t=2}(C_{27}(1,3,8,10))$ is displayed as shown in Figure 8. Moreover, it indicates  $C_{27}(2,3,7,11)$ = $\theta_{27,3,2}(C_{27}(1,3,8,10))$ is Type-2 isomorphic w.r.t. $m$ = 3 (to $C_{27}(1,3,8,10)$).

\item [\rm (v)] In order to get Adam's isomorphic circulant graph of $C_{27}(1,3,8,10)$, after step (iii), choose AdamIso and click the button, circulant graph $C_{27}(2,6,7,11)$ is displayed as shown in Figure 9. It indicates  $C_{27}(2,6,7,11)$ = $C_{27}(2(1,3,8,10))$ is Type-1 isomorphic to $C_{27}(1,3,8,10)$. 

When we click the Continue button, circulant graph $C_{27}(4,5,12,13)$ is displayed as shown in Figure 10. It indicates  $C_{27}(4,5,12,13)$ = $C_{27}(4(1,3,8,10))$ is Type-1 isomorphic to $C_{27}(1,3,8,10)$. 

Repeat this step to get more such graphs, if exist.

\item [\rm (vi)] Click Exist when you want to come out of the trail.

\item [\rm (vii)] Repeat the above steps for another circulant graph. And so on.
\end{enumerate}
 
%Fig 7-4
%%%%%%%%%%%%%%%%%%%%%%%%%%%
\begin{figure}[ht]
\centerline{\includegraphics[width=4.1in]{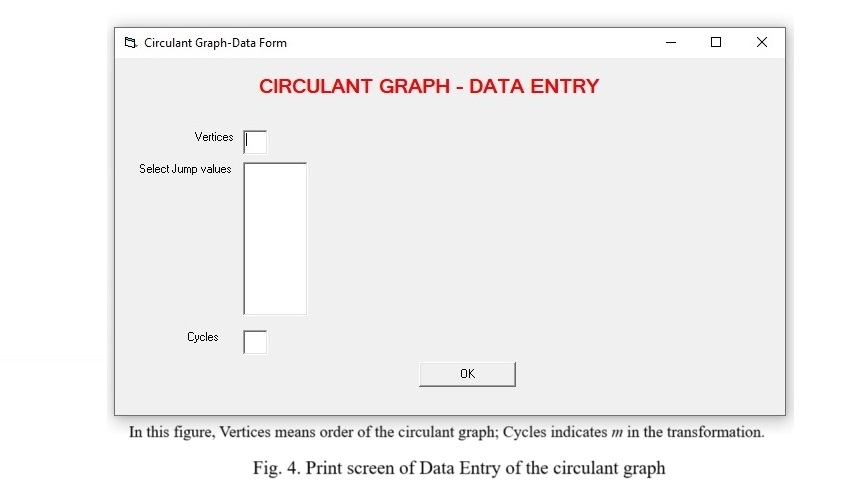}}
\end{figure}
%%%%%%%%%%%%%%%%%%%%%%%%%%%%%%%%%
%Fig 7-5
%%%%%%%%%%%%%%%%%%%%%%%%%%%
\begin{figure}[ht]
\centerline{\includegraphics[width=3.9in]{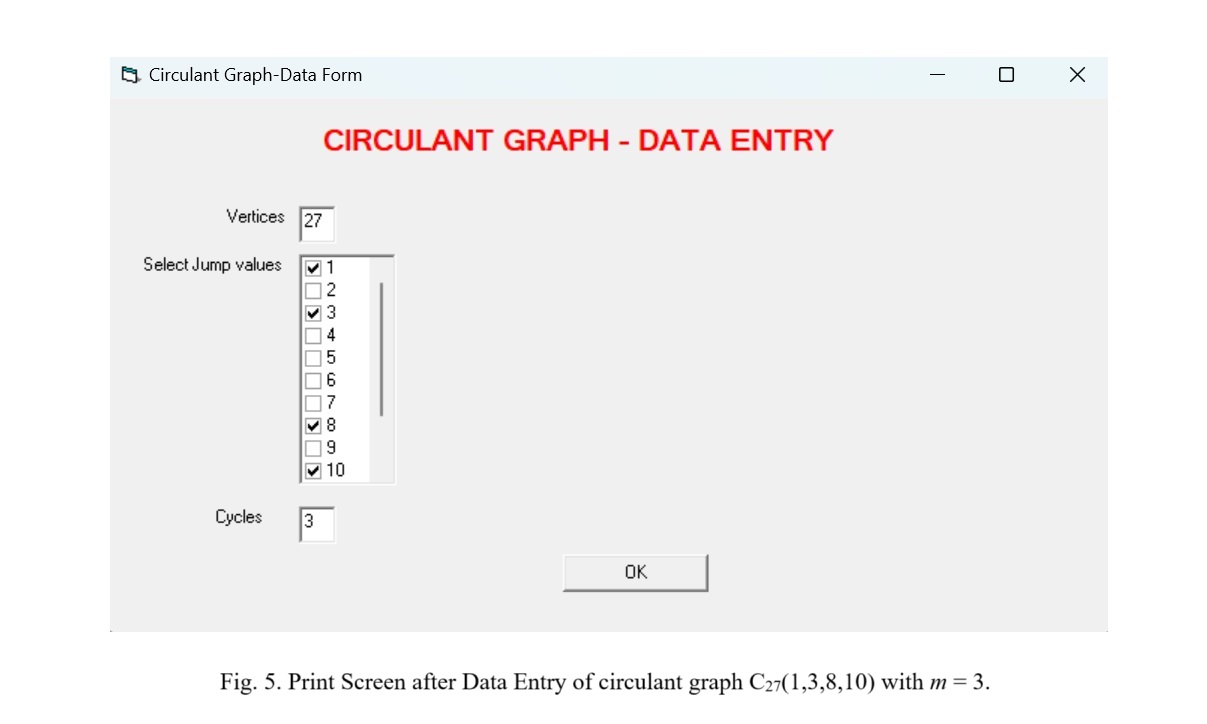}}
\end{figure}
%%%%%%%%%%%%%%%%%%%%%%%%%%%%%%%%%

%Fig 7-6
%%%%%%%%%%%%%%%%%%%%%%%%%%%
\begin{figure}[ht]
\centerline{\includegraphics[width=3.9in]{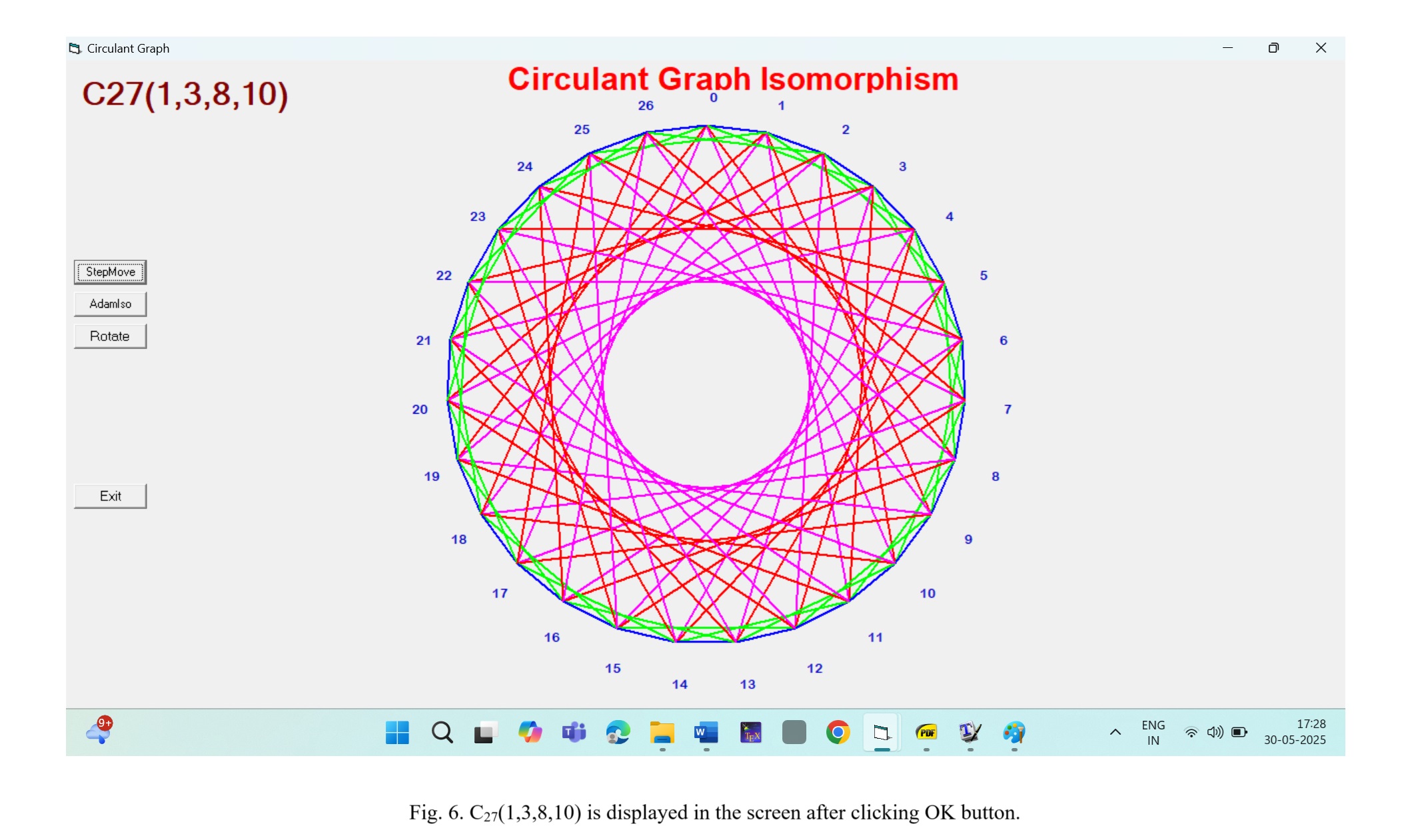}}
\end{figure}
%%%%%%%%%%%%%%%%%%%%%%%%%%%%%%%%%

%Fig 7-7
%%%%%%%%%%%%%%%%%%%%%%%%%%%
\begin{figure}[ht]
\centerline{\includegraphics[width=4.1in]{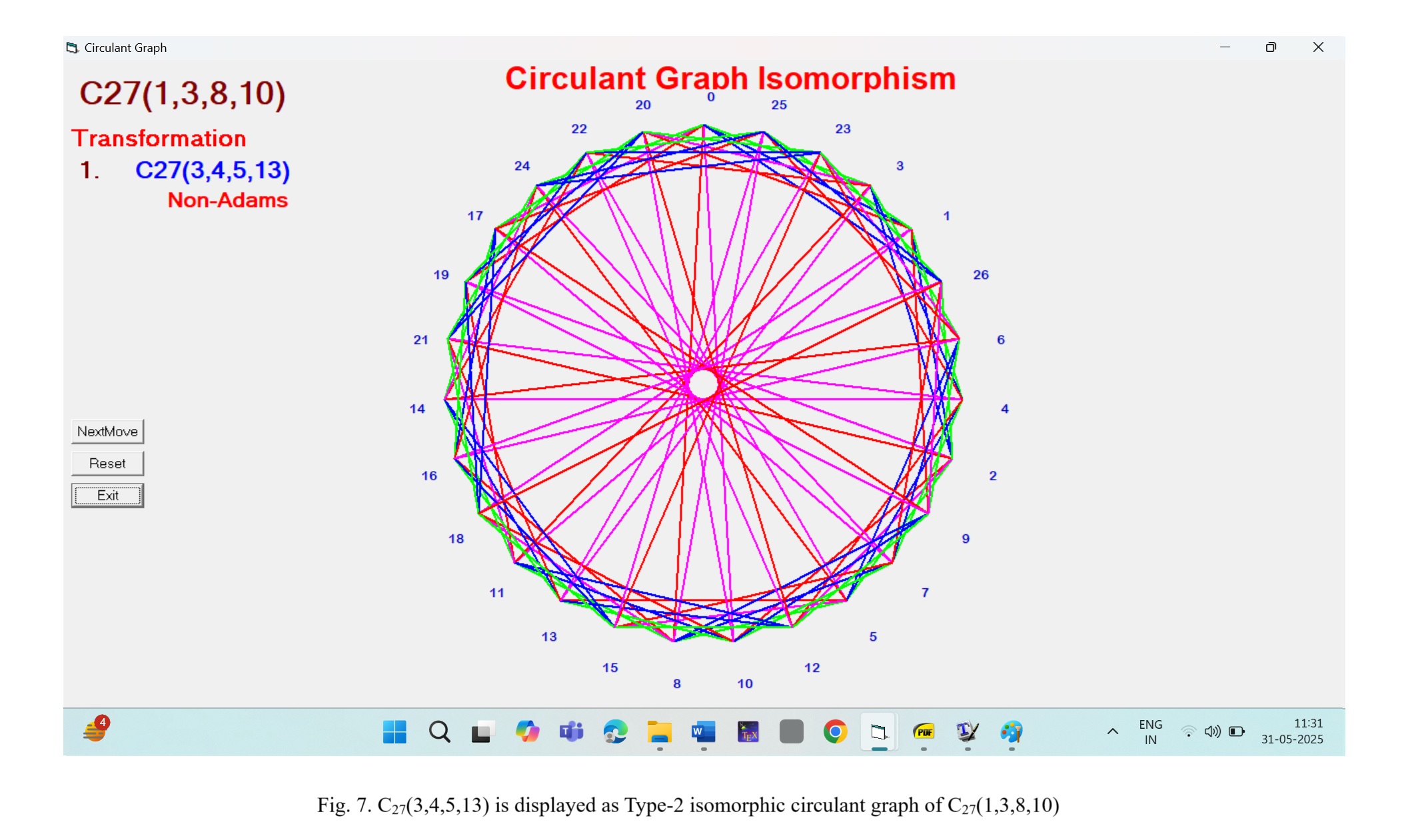}}
\end{figure}
%%%%%%%%%%%%%%%%%%%%%%%%%%%%%%%%%

%Fig 7-8
%%%%%%%%%%%%%%%%%%%%%%%%%%%
\begin{figure}[ht]
\centerline{\includegraphics[width=4.04in]{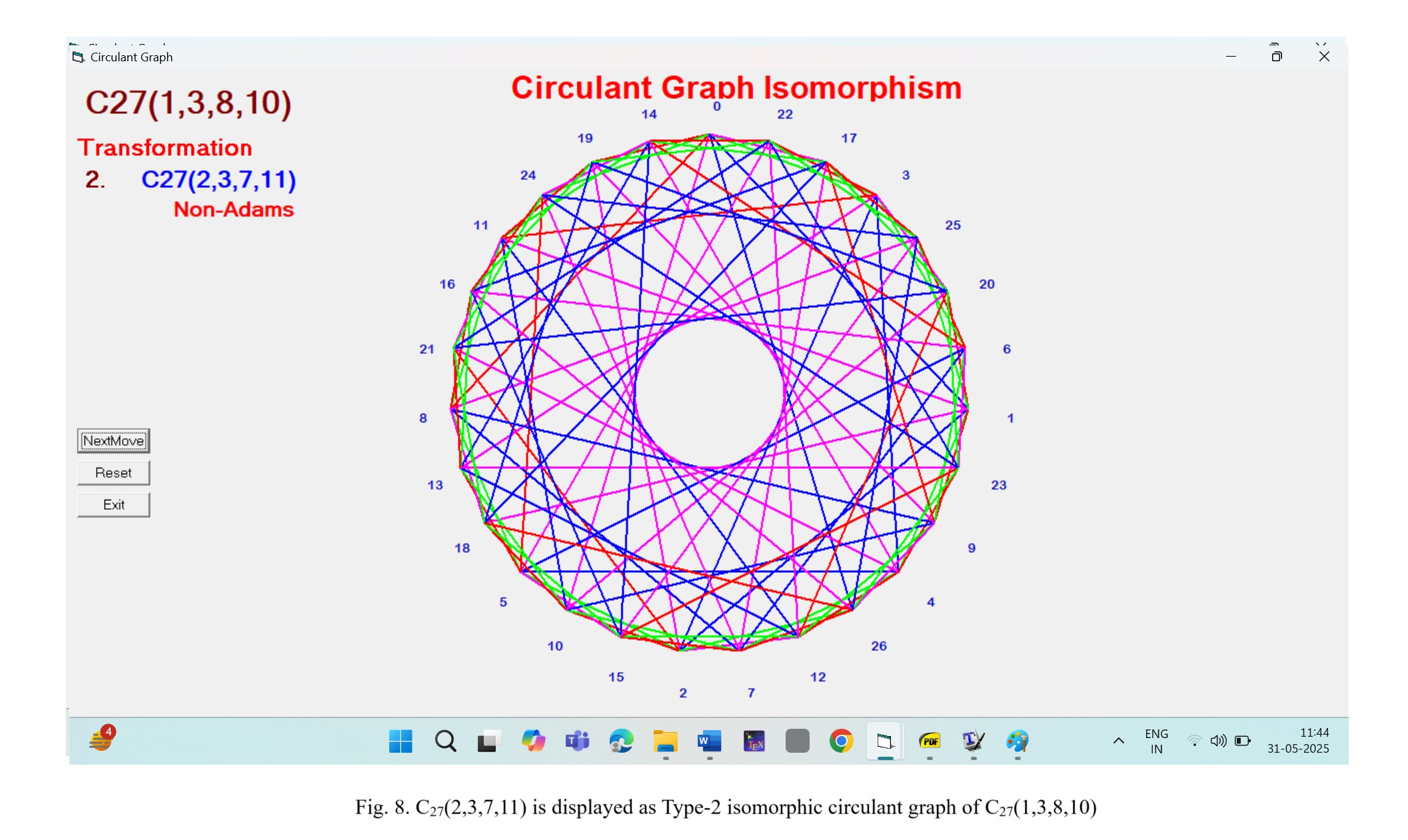}}
\end{figure}
%%%%%%%%%%%%%%%%%%%%%%%%%%%%%%%%%

%Fig 7-9
%%%%%%%%%%%%%%%%%%%%%%%%%%%
\begin{figure}[ht]
\centerline{\includegraphics[width=4.08in]{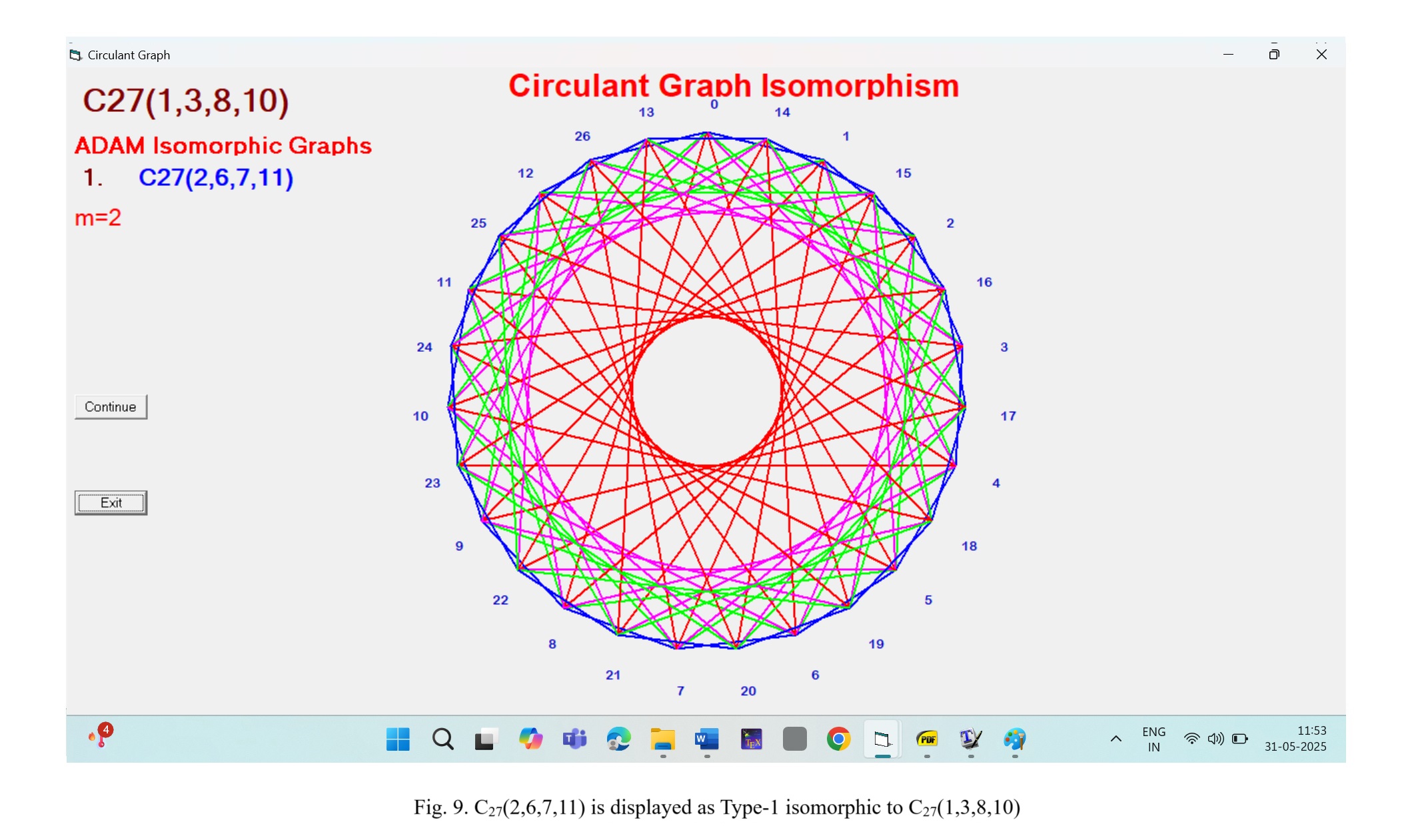}}
\end{figure}
%%%%%%%%%%%%%%%%%%%%%%%%%%%%%%%%%

%Fig 7-10
%%%%%%%%%%%%%%%%%%%%%%%%%%%
\begin{figure}[ht]
\centerline{\includegraphics[width=4.3in]{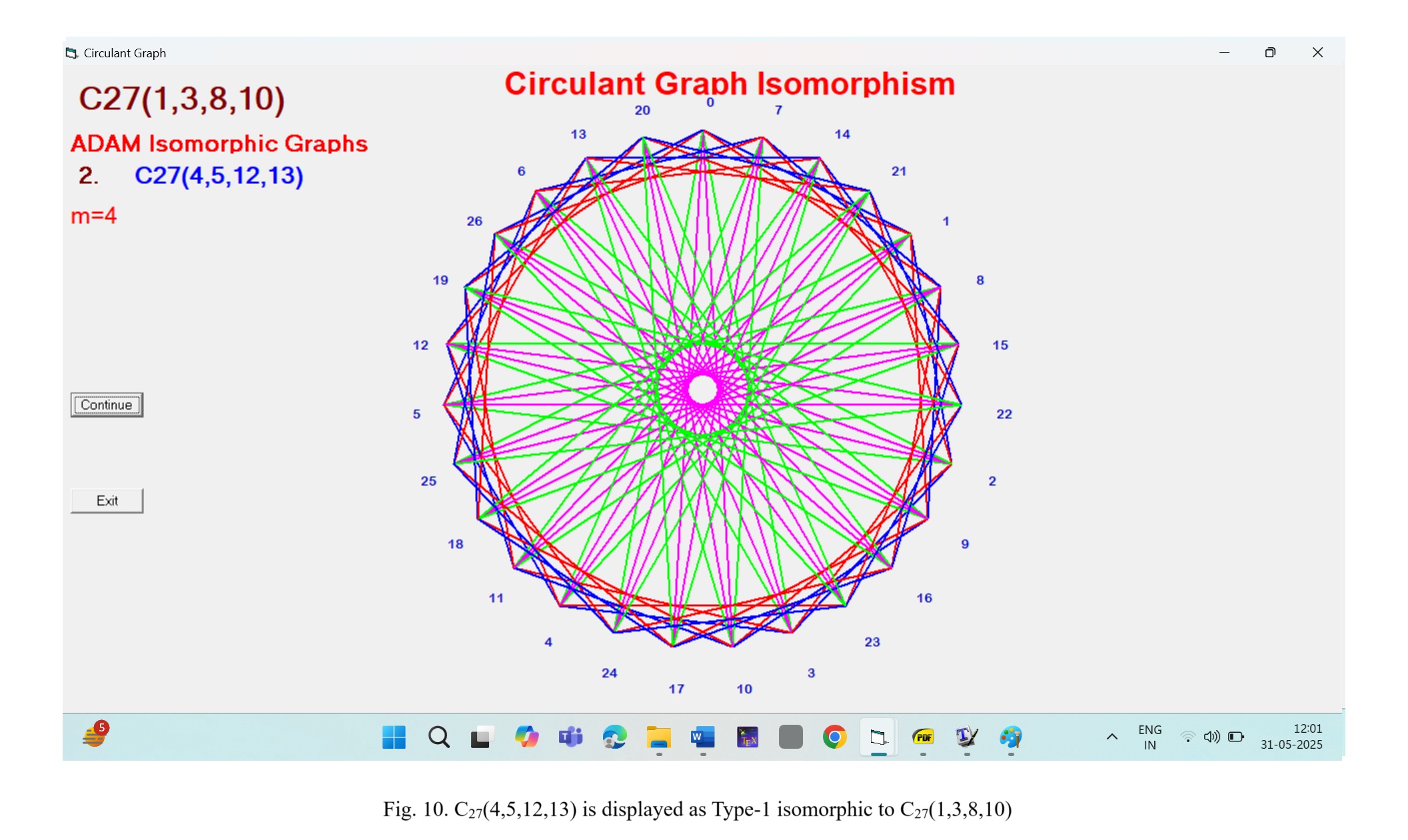}}
\end{figure}
%%%%%%%%%%%%%%%%%%%%%%%%%%%%%%%%%

\section {Conclusion}
In general, computational complexity on checking isomorphism of graphs increases with order and degree of the graphs when the graphs are regular, except complete graphs. When the graphs are not regular, then their computational complexity in checking isomorphism can be reduced by considering their disjoint maximal induced regular subgraphs which partition their vertex sets and their interconnected subgraphs. The authors feel that computer programs similar to POLY415.EXE can be used to attack problems on graph isomorphism and on crossing numbers. Also, from different isomorphic circulant graphs of a given circulant graphs, it is possible to learn more on edge coloring, cycle coloring, edge decomposition and crossing number of circulant graphs.   

\vspace{.3cm}
\noindent
\textbf{Conflict of Interest}

\noindent
\textit {The authors declare that there is no conflict of interests regarding the publication of this paper.}

\vspace{.3cm}
\noindent
{\bf Acknowledgement.}\quad {\rm We express our sincere thanks Dr. A. Christopher and Mr. R. Satheesh of S.T. Hindu College, Nagercoil, Tamil Nadu, India for their assistance to develop
the VB programs. We also express our gratitude to the Central University of Kerala, Kasaragod, Kerala; St. Jude’s College, Thoothoor and S. T. Hindu College, Nagercoil; and Lerroy Wilson Foundation, India (www.WillFoundation.co.in) for providing facilities to do this research work.}

\end{document}